\documentclass[12pt]{article}

\usepackage{graphicx,psfrag}
\usepackage{psfrag,eepic,epsfig}
\usepackage{sectsty}
\usepackage{setspace}
\usepackage{color}

\usepackage{setspace}
\usepackage{amsmath, epsfig, cite, ifpdf}
\usepackage{amssymb}
\usepackage{amsfonts}
\usepackage{latexsym}
\usepackage{graphicx,psfrag,eepic,rawfonts}
\usepackage{enumerate}
\usepackage{indentfirst}
\input{epsf}

\newtheorem{thm}{Theorem}[section]

\newtheorem{defi}[thm]{Definition}
\newtheorem{lem}[thm]{Lemma}
\newtheorem{conj}[thm]{Conjecture}

\newcommand{\qed}{{\hfill\rule{4pt}{7pt}}}
\def\pf{\noindent {\it Proof.} }

\numberwithin{equation}{section}

\makeatletter \@addtoreset{equation}{section} \makeatother

\title {\bf The matching energy of random graphs\footnote{Supported
by NSFC No.11371205, the ``973" program No.2013CB834204, and PCSIRT.}}

\author{
{\small Xiaolin Chen, Xueliang Li}\\
{\small Center for Combinatorics and LPMC-TJKLC}\\
{\small Nankai University, Tianjin 300071, China}\\
{\small E-mail: chxlnk@163.com; lxl@nankai.edu.cn}\\
{\small Huishu Lian}\\
{\small College of Sciences}\\
{\small China University of Mining and Technology, Xuzhou 221116, China}\\
{\small E-mail: lhs6803@126.com}   }
\date{}
\begin{document}

\maketitle

\begin{abstract}
The matching energy of a graph was introduced by Gutman and Wagner,
which is defined as the sum of the absolute values of the roots of
the matching polynomial of the graph. For the random graph $G_{n,p}$
of order $n$ with fixed probability $p\in (0,1)$, Gutman and Wagner
[I. Gutman, S. Wagner, The matching energy of a graph, Discrete
Appl. Math. 160(2012), 2177--2187] proposed a conjecture that the
matching energy of $G_{n,p}$ converges to
$\frac{8\sqrt{p}}{3\pi}n^{\frac{3}{2}}$ almost surely. In this
paper, using analysis method, we prove that the conjecture is true.
\\


\end{abstract}

\section{Introduction}

Let $G$ be a finite simple graph of order $n$ with vertex set $V(G)$
and edge set $E(G)$. A matching of $G$ is a set of independent edges
in $G$, and an $r$-matching of $G$ is a matching of $G$ that has
exactly $r$ edges. By $m_r(G)$ we denote the number of $r$-matchings
in $G$. It is easy to verify that for $k<0$ and $k>\lceil
n/2\rceil$, $m_r(G)=0$. And when $r=1$, $m_1(G)$ is the size of $G$.
For convenience, we define $m_{0}(G)=1$. The matching polynomial
$m(G,x)$ \cite{FAR,GG,GUT2} of a graph $G$ is defined as
\begin{equation*}
m(G,x)=\sum\limits_{k=0}^{\lfloor n/2\rfloor}(-1)^{k}m_k(G)x^{n-2k}.
\end{equation*}

The matching polynomial has been widely studied and many results on
the properties of the matching roots have been obtained; see
\cite{FAR,GOD1,GOD2,GUT2,GG,HL}. For any graph $G$, all the matching
roots are real. If $\lambda$ is a matching root, then $-\lambda$ is
also a matching root. That is, the matching roots are symmetric.
Moreover, the matching polynomial has many important implications in
statistical physics and chemistry; see \cite{GUT1,GMT,HL}.

In \cite{GW}, Gutman and Wagner introduced the matching energy (ME)
of a graph $G$, which is defined as the sum of the absolute values
of the roots of the matching polynomial of $G$. Note that the
concept of the energy $E(G)$ of a simple undirected graph $G$ was
introduced by Gutman in \cite{GUT}. Afterwards, there have been lots
of research papers on this topic. A systematic study of this topic
can be found in the book \cite{LSG}. In \cite{GW}, Gutman and Wagner
pointed out that the matching energy is a quantity of relevance for
chemical applications. Moreover, they arrived at the simple relation
$$TRE(G)=E(G)-ME(G),$$
where $TRE(G)$ is the so-called ``topological resonance energy" of
the graph $G$. For more information about the applications of the
matching energy, we refer the reader to \cite{GUT1,GMT}. Recently,
there have been some results on the extremal values of the matching
energy of graphs; see \cite{CS,JLS,LZS,LY}.

When we add an edge to a graph, the matching energy increases
strictly.
\begin{lem}(\cite{GW})\label{add}
Let $G$ be a graph and $e$ one of its edges. Let $G-e$ be the
subgraph obtained by deleting from $G$ the edge $e$, but keeping all
the vertices of $G$. Then
\begin{equation*}
ME(G-e)<ME(G).
\end{equation*}
\end{lem}
Therefore, the complete graph $K_n$ attains the maximum matching energy among all graphs of order $n$.
In \cite{GW}, Gutman and Wagner got the following lemma which gives an asymptotic estimation of
the matching energy of $K_n$.
\begin{lem}(\cite{GW})
The matching energy of the complete graph $K_n$ is asymptotically equal to $\frac{8}{3\pi}n^{3/2}$.
More precisely,
\begin{equation}
ME(K_n)= \frac{8}{3\pi}n^{3/2}+O(n).
\end{equation}
\end{lem}

The above lemma can been thought as the upper bound of the matching energy of all graphs of order $n$.
Moreover, they studied the the lower bound of the matching energy of random graphs. Now we recall some notion in probability, we say an event holds almost surely (a.s.) if it occurs with probability 1. An event holds asymptotically almost surely (a.a.s.) if the probability of success goes to 1 as $n\rightarrow \infty$.
\begin{lem}(\cite{GW})
Consider the random graph $G_{n,p}$ of order $n$ with fixed probability $p\in (0,1)$. Then
\begin{equation}
ME(G_{n,p})\geq \frac{\sqrt{p}}{\pi}n^{3/2}+O(\sqrt{n}\ln n)
\end{equation}
holds asymptotically almost surely.
\end{lem}

Based on the above analysis, they conjectured that
\begin{conj}(\cite{GW})\label{conj}
For any fixed probability $p\in (0,1)$,
$$n^{-3/2}ME(G_{n,p})\longrightarrow \frac{8\sqrt{p}}{3\pi}$$
asymptotically almost surely.
\end{conj}

This paper is to confirm the conjecture. The rest of the paper is
organized as follows. In Section 2, we introduce the empirical
matching distribution and list our main results: the empirical
matching distribution converges weakly to the semicircle
distribution; the asymptotic formula of the matching energy of
random graphs. The explicit proofs will be shown in Sections
\ref{PTH1} and \ref{PTH2}. Throughout the paper we use the following
standard asymptotic notation: as $n\rightarrow \infty$,
$f(n)=o(g(n))$ means that $f(n)/g(n)\rightarrow 0$; $f(n)=O(g(n))$
means that there exists a constant $C$ such that $|f(n)|\leq C
g(n)$.

\section{Matching energy of random graphs}

In this section, we present our main results of this paper.
\begin{defi}
For the random graph $G_{n,p}$ of order $n$ with fixed probability
$p\in (0,1)$, let $x_1(G_{n,p})\geq \cdots\geq x_n(G_{n,p})$ be the
roots of the matching polynomial $m(G_{n,p},x)$, since all roots of
the matching polynomial are real. Then let $\lambda_i(G_{n,p})=
\frac{1}{\sqrt{np}}x_i(G_{n,p})$ for all $1\leq i\leq n$.

We define the empirical matching distribution (EMD) as a
distribution function $F_{n}(x)$ where
\begin{equation*}
F_{n}(x)=\frac{1}{n}\Big|\{\lambda_i(G_{n,p})|\lambda_i(G_{n,p})\leq x,i=1,2,\ldots,n\}\Big|
\end{equation*}
\end{defi}

The empirical matching distribution can be thought as the root
distribution of the matching polynomial. Most work on the root
distribution focuses on the spectral distributions of random
matrices. The study can be traced back to the pioneer work
semicircle law discovered by Wigner in \cite{Wigner}. Afterwards,
the research about the spectral distributions of many sorts of
random matrices became the topics in mathematics and physics. For
more details, we refer the reader to books \cite{AGZ,BS,TT}. In this
paper, we find that the empirical matching distribution has the
similar convergent property.
\begin{thm}\label{THM1}
For the random graph $G_{n,p}$ of order $n$ with $p\in (0,1)$,
the empirical matching distribution
$F_{n}(x)$ converges weakly almost surely to the standard semicircle  distribution $F_{sc}(x)$,
whose density is given by
\begin{equation*}
\rho_{sc}(x):=\frac{1}{2\pi}\sqrt{4-x^2}\,\,\,_{|x|\leq 2}.
\end{equation*}
That is, for any bounded continuous function $f$ on $R$,
\begin{equation}
\int_{R}fdF_{n}(x) \longrightarrow \int_{R}fdF_{sc}(x)
\end{equation}
almost surely.
\end{thm}

From the above theorem, we can estimate the matching energy of the
random graph $G_{n,p}$. Before proceeding, we should note that the
semicircle law has been used to study many energies, such as the
energy in \cite{NIK, DLL3}, the Laplacian energy in \cite{DLL1}, the
skew energy in \cite{CLL}, and other energies in \cite{DLL2}.

We prove Conjecture \ref{conj} by the following theorem.
\begin{thm}\label{THM2}
For $p\in (0,1)$, the matching energy $ME(G_{n,p})$ of the random
graph $G_{n,p}$ enjoys asymptotically almost surely the following
equation:
\begin{equation*}
ME(G_{n,p})=n^{3/2}p^{1/2}\left(\frac{8}{3\pi}+o(1)\right).
\end{equation*}
\end{thm}

\section{Empirical matching distributions (EMDs) of
random graphs}\label{PTH1}

In order to prove the empirical matching distribution (EMD) of
$G_{n,p}$ converges weakly to the standard semicircle distribution,
we utilize the so-called moment method. The moment method has been
used extensively in random matrices, specially in semicircle
distribution \cite{AGZ,BS,TT}. The key point of moment method is to
show that the moments of EMDs converge almost surely to the moments
of the semicircle law. Thus, in this section, we just need to verify
the almost sure convergence of the moments of the matching roots.
For the gap between Theorem \ref{THM1} and the moment method, we can
fill it by the method similar to random matrices (\cite{AGZ},P.11).
Thus, we do not show it in this section. However, for convenience of
the reader, we present it in Appendix \ref{Appendix}.

In Subsection \ref{1}, we review the relationship between the
tree-like closed walks and the moments of matching roots. In
Subsection \ref{2}, we prove a weaker convergence named convergence
in expectation. In Subsection \ref{3}, using some probabilistic
inequalities, we upgrade the convergence in expectation to the
almost sure convergence.

\subsection{Tree-like walks of graphs}\label{1}

In \cite{GOD2}, Godsil introduced a new type of closed walks, namely
the tree-like walks, which have very closed relation to the matching
roots. Before proceeding, we recall some notation. A closed walk is
called minimal if only the first vertex and the last vertex
coincide. For any closed walk $w$ with length nonzero, we can
uniquely decompose it into the form $\alpha\beta\gamma$, where
$\alpha$ is a path, $\beta$ is a minimal closed walk and has no
common vertex with $\alpha$. Then, we call $\beta$ the first minimal
closed walk in $w$. Suppose the first vertex in $\beta$ is $u_0$. By
deleting the closed walk $\beta$, we get a new closed walk
$w'=\alpha u_0\gamma$. We may again decompose the new closed walk
$w'$ and get the first minimal closed walk in $w'$. By continuing
the operation in this way, we can get a sequence of minimal closed
walks. And we call the members of this sequence the factors of the
closed walk $w$.

For example, let $w$ be a closed walk $\{u,a,b,c,d,c,d,u\}$, the
first closed walk of $w$ is $\{c,d,c\}$. And by deleting it, we get
a new closed walk $w'=\{u,a,b,c,d,u\}$ which is a minimal closed
walk. Then $\{c,d,c\}$ and $\{u,a,b,c,d,u\}$ are two factors of $w$.

When the closed walk $w$ is in a tree, we find that all the factors
of $w$ have length two. A closed walk is called tree-like if all its
factors have length two. In \cite{GOD2}, Godsil gave a combinatorial
interpretation about the moments of matching roots.
\begin{lem}(\cite{GOD2})\label{Moment}
Let $m(G,x)$ be the matching polynomial of a graph $G$. Then the
rational function $xm'(G,x)/m(G,x)$ is the generating function, in
the variable $x^{-1}$, for the closed tree-like walks in $G$.
\end{lem}
By Lemma \ref{Moment}, we can get the following equation.
\begin{equation}
\frac{xm'(G,x)}{m(G,x)}=\sum\limits_{i=1}^{n}\frac{x}{x-x_i}
=\sum\limits_{k=0}^{\infty}(\sum_{i=1}^{n}x_i^k)x^{-k}
\end{equation}
Then $\sum\limits_{i=1}^{n}x_i^k$ is the number of closed tree-like
walks with length $k$ in $G$.

\subsection{Convergence of the moments in expectation}\label{2}

We recall a weaker notion of convergence, named convergence in
expectation, which is defined as follows. We say that EMDs converge
in expectation to the semicircle law, if
\begin{equation}
\mathbf{E}\int_{R}f(x)dF_{n}(x)\rightarrow\int_{R}f(x)dF_{sc}(x)\,\,\,a.s.
\end{equation}
as $n\rightarrow \infty$,
for all bounded continuous function $f(x)$.
\begin{thm}\label{MonCon}
For any positive integer $k$, $\lim\limits_{n\rightarrow \infty}
\mathbf{E}\int_{R}x^kdF_{n}(x)=\int_{R}x^kdF_{sc}(x)\,\,\,a.s.$
\end{thm}
\pf For any positive integer $k$, the expected $k$-th moment of the EMD is
\begin{equation}
\mathbf{E}\int_R x^kdF_{n}(x)=\mathbf{E}\frac{1}{n}\sum\limits_{i=1}^{n}\lambda_i^k,
\end{equation}
and the $k$-th moment of the standard semicircle distribution is
\begin{equation*}
\int_{-2}^{2}x^kdF_{sc}(x).
\end{equation*}
Next, we just need to prove that for every fixed integer $k$,
\begin{equation*}
\mathbf{E}\frac{1}{n}\sum\limits_{i=1}^{n}\lambda_i^k \longrightarrow
\int_{-2}^{2}x^kdF_{sc}(x), \text{ as } n\rightarrow \infty.
\end{equation*}
And by Lemma \ref{Moment}, we have that
$$\mathbf{E}\frac{1}{n}\sum\limits_{i=1}^{n}\lambda_i^k=
\frac{1}{n^{1+k/2}p^{k/2}}\sum\limits_{w}\mathbf{E}(X(w)),$$ where
$w$ denotes a tree-like closed walk with length $k$ in the random
graph $G_{n,p}$, and $X(w)$ is a random variable taking value $1$ if
$w$ occurs and $0$ otherwise.

When $k$ is an odd number, then $\int_{-2}^{2}x^kdF_{sc}(x)=0$.
Since $w$ is a tree-like closed walk, the length of $w$ must be even
and for any edge of $G_{n,p}$, the total number of times that this
edge appears in $w$ is even. Then
$\frac{1}{n^{1+k/2}p^{k/2}}\sum\limits_{w}\mathbf{E}(X(w))=0$.

Now considering the even number $k=2m$ where $m\geq 1$. We have
\begin{align*}
\int_{-2}^{2}x^kdF_{sc}(x)&=\frac{1}{2\pi}\int_{-2}^{2}x^{k}\sqrt{4-x^2}dx=
\frac{1}{\pi}\int_{0}^{2}x^{2m}\sqrt{4-x^2}dx\\
&=\frac{2^{2m+1}}{\pi}\cdot\frac{\Gamma(m+1/2)\Gamma(3/2)}{\Gamma(m+2)}=\frac{1}{m+1}{2m
\choose m},
\end{align*}
where $\Gamma(x)$ is the standard Gamma function. Let $t$ denote the
number of distinct vertices in a tree-like walk $w$. It is easy to
see that $t$ is no more than $m+1$. We divide the discussion into
two cases.

\textbf{Case 1.} For a tree-like walk $w$ with $t\leq m$, it is easy
to check that the number of this kind of closed walks is less than
$t^k$ and the edges of $w$ induce a connected graph $H_w$ of order
$t$ in $G_{n,p}$. Then we have
\begin{align*}
&\frac{1}{n^{1+k/2}p^{k/2}}\sum_{t=1}^{m}\sum_{|w=\{i_1,\ldots,i_k\}|=t}\mathbf{E}(X(w))\\
\leq&\frac{1}{n^{1+m}p^m}\sum_{t=1}^{m}n^t\cdot t^k\cdot \mathbf{E}(X(H_w))\\
\leq&\frac{1}{n^{1+m}p^m}\sum_{t=1}^{m}n^t\cdot t^k\cdot p^{t-1}\\
\leq&\frac{1}{n^{1+m}p^m}\cdot m \cdot n^m\cdot m^{2m}\cdot p^{m-1}\\
=&\frac{m^{2m+1}}{np}=O\left(\frac{1}{np}\right).
\end{align*}

\textbf{Case 2.} For a tree-like walk $w$ with $t=m+1$, we know that
each edge of the closed walk $w$ appears twice, and the number of
distinct edges in the closed walk $w$ is $m$. The number of such
kind of closed walks can be determined by the following lemma.
\begin{lem}\cite{BS}
The number of the closed walks of length $2m$ which satisfy that
each edge and its inverse edge both appear once in the closed walks
is $\frac{1}{m+1}{2m \choose m}$.
\end{lem}
Then, we obtain
\begin{align*}
&\frac{1}{n^{1+k/2}p^{k/2}}\sum_{|w=\{i_1,\ldots,i_k\}|=m+1}\mathbf{E}(X(w))\\
\leq&\frac{1}{n^{1+k/2}p^{k/2}}{n \choose m+1}\frac{1}{m+1}{2m \choose m}p^{m}\\
=&(1-\frac{1}{n})(1-\frac{2}{n})\cdots(1-\frac{m}{n})\frac{1}{m+1}{2m \choose m}.
\end{align*}
When $n$ is large enough, we have
$$\frac{1}{n^{1+k/2}p^{k/2}}\sum_{|w=\{i_1,\ldots,i_k\}|=m+1}\mathbf{E}(X(w))=(1+o(1))\frac{1}{m+1}{2m \choose m}.$$
By the above analysis, we have
\begin{equation*}
\frac{1}{n^{1+k/2}p^{k/2}}\sum\limits_{w}\mathbf{E}(X(w))=
\begin{cases}
0 & \text{ if }k=2m+1;\\ \frac{1}{m+1}{2m \choose m}(1+o(1))+O\left(\frac{1}{np}\right) & \text{ if }k=2m,
\end{cases}
\end{equation*}
Since $np\rightarrow \infty$ as $n\rightarrow \infty$, it follows that when $n\rightarrow \infty$,
\begin{equation*}
\mathbf{E}\int_{R}x^kdF_{n}(x)\longrightarrow\int_{R}x^kdF_{sc}(x)\,\,\,a.s.
\end{equation*}
The proof is thus complete. \qed

\subsection{Almost sure convergence of the moments}\label{3}

In this section, we prove the following theorem.
\begin{thm}\label{THM0}
For any positive integer $k$, $\lim\limits_{n\rightarrow \infty}
\int_{R}x^kdF_{n}(x)=\int_{R}x^kdF_{sc}(x)\,\,\,a.s.$
\end{thm}
\pf To prove that the moments of EMDs converge almost surely to the
moments of the semicircle distribution, by Borel-Cantelli Lemma
(\cite{B},P.60), it will be sufficient to show that
\begin{equation}\label{Converge}
\sum\limits_{n=1}^{\infty}P\left(|\int_R x^kdF_n(x)-\int_R x^k dF_{sc}(x)|>\epsilon\right)<\infty.
\end{equation}
By the triangle inequality, we have
\begin{align*}
&P\left(|\int_R x^kdF_n(x)-\int_R x^k dF_{sc}(x)|>\epsilon\right)\\
\leq&P\left(|\int_R x^kdF_n(x)-\mathbf{E}\int_R x^kdF_n(x)|+|\mathbf{E}\int_R x^kdF_n(x)-\int_R x^k dF_{sc}(x)|>\epsilon\right)\\
\leq&P\left(|\int_R x^kdF_n(x)-\mathbf{E}\int_R x^kdF_n(x)|>\epsilon/2\right)+P\left(|\mathbf{E}\int_R x^kdF_n(x)-\int_R x^k dF_{sc}(x)|>\epsilon/2\right)
\end{align*}
By Theorem \ref{MonCon}, we have
$$\sum\limits_{n=1}^{\infty} P\left(|\mathbf{E}\int_R x^kdF_n(x)-\int_R x^k dF_{sc}(x)|>\epsilon/2\right)<\infty.$$
By Chebyshev's inequality, we have
\begin{equation*}
P\left(|\int_R x^kdF_n(x)-\mathbf{E}\int_R x^kdF_n(x)|>\epsilon/2\right)<\frac{4}{\epsilon^2}\mathbf{Var}\int_R x^k dF_n(x).
\end{equation*}
Next, we consider the variance of the moments of EMDs. Let $E(w)$ be the edge set of $w$ and $V(w)$ the vertex set of $w$. Then
\begin{align*}
\mathbf{Var}\int_R x^k dF_n(x)=\,\,\,& \mathbf{E}(\int_R x^k dF_n(x))^2-(\mathbf{E}\int_R x^k dF_n(x))^2\\
=\,\,\,&\frac{1}{n^{2+k}p^k}\sum\limits_{w_1,w_2}
\left(\mathbf{E}(X(w_1, w_2))-\mathbf{E}(X(w_1))\mathbf{E}(X(w_2))\right),
\end{align*}
where $w_1,w_2$ are both tree-like closed walks, and $X(w_1,w_2)$ is a random variable taking value $1$ if $w_1,w_2$ both occur and $0$ otherwise. If the tree-like closed walks $w_1$ and $w_2$ are edge-disjoint, then $\mathbf{E}(X(w_1, w_2))=\mathbf{E}(X(w_1))\mathbf{E}(X(w_2))$.
Hence, we only need to consider the pairs of tree-like closed walks which share at least one edge.
Then $|V(w_1)\cup V(w_2)|\leq |V(w_1)|+|V(w_2)|-2$. Since $|V(w)|\leq k/2+1$, the number of pairs of tree-like closed walks which contribute to the sum is no more than $n^k$.
It follows that
\begin{equation*}
\frac{1}{n^{2+k}p^k}\sum\limits_{w_1,w_2}(\mathbf{E}(X(w_1, w_2))-\mathbf{E}(X(w_1))\mathbf{E}(X(w_2)))
=\frac{1}{n^{2+k}}n^{k}O(1)=O(n^{-2}).
\end{equation*}
and
\begin{equation*}
P\left(|\int_R x^kdF_n(x)-\mathbf{E}\int_R x^kdF_n(x)|>\epsilon/2\right)=\frac{4}{\epsilon^2}O(n^{-2}).
\end{equation*}
Then for each $\epsilon>0$, we have \,\,\,$\sum\limits_{n=1}^{\infty}\frac{4}{\epsilon^2}O(n^{-2})<\infty.$

Thus, Equation \ref{Converge} follows and the proof is complete. \qed

\section{The proof of Theorem \ref{THM2}}\label{PTH2}

In this section, we give a proof of Theorem \ref{THM2}. Before proceeding, we recall a useful lemma.
\begin{lem}(\cite{B},P.198)\label{UNI}
If the distribution function $F_n$  converges weakly to the distribution function $F$,
and $F$ is everywhere continuous, then $F_n(x)$ converges to $F(x)$ uniformly in $x$.
\end{lem}
\textbf{Proof of Theorem \ref{THM2}:} For positive constants $M$ and $\delta$, we define a bounded continuous function
\begin{equation}
f(x)=
\begin{cases}
|x| & \text{ if } -M\leq x\leq M;\\
M-\frac{M}{\delta}(x-M) & \text{ if }M<x\leq M+\delta;\\
 M+\frac{M}{\delta}(x+M) & \text{ if }-M-\delta\leq x<-M;\\
0 & \text{otherwise.}
\end{cases}
\end{equation}
Then by Theorem \ref{THM1}, we get
\begin{equation}
\lim\limits_{n\rightarrow \infty}\int_{R}f(x)dF_{n}(x) =\int_{R}f(x)dF_{sc}(x)\,\,\,\,a.s.
\end{equation}
Let $M>4$ and $\delta>0$. And it follows that $\int_{|x|\geq M}f(x)dF_{sc}(x)=0$.
Now we consider the following inequalities.
\begin{align*}
&\big|\int\limits_{|x|\leq M}|x|dF_{n}(x)-\int\limits_{|x|\leq M}|x|dF_{sc}(x)\big|\\
&\leq\big|\int_R f(x)dF_{n}(x)-\int_R f(x)dF_{sc}(x)\big|+2\big|\int\limits_{M<x\leq M+\delta} f(x)dF_{n}(x)\big|
\end{align*}
where
\begin{align*}
&2\big|\int\limits_{M<x\leq M+\delta} f(x)dF_{n}(x)\big|\leq2M\big|\int\limits_{M<x\leq M+\delta} dF_{n}(x)\big|\\
&\leq2M\big|F_{n}(M+\delta)-F_{sc}(M+\delta)+F_{sc}(M)-F_{n}(M)+\int\limits_{M<x\leq M+\delta} dF_{sc}(x)\big|\\
&\leq2M\Big(\big|F_{n}(M+\delta)-F_{sc}(M+\delta)\big|+\big|F_{sc}(M)-F_{n}(M)\big|+\big|\int\limits_{M<x\leq
M+\delta} dF_{sc}(x)\big|\Big).
\end{align*}

From Equation (2.1), we know that for any small number
$\varepsilon>0$, there exists a number $N_1$ such that $\big|\int_R
f(x)dF_{n}(x)-\int_R f(x)dF_{sc}(x)\big|<\frac{\varepsilon}{3}$ for
all $n>N_1$. Since $F_{sc}(x)$ is everywhere continuous, by Lemma
\ref{UNI}, there exist a number $N_2>0$ and a number $\delta_1>0$
such that, when $n>N_2$ and $0<\delta<\delta_1$,
$\big|F_{n}(x)-F_{sc}(x)\big|<\frac{\varepsilon}{12M}$ for all $x$,
and  $\big|\int\limits_{M<x\leq M+\delta}
dF_{sc}(x)\big|<\frac{\varepsilon}{12M}$. Hence, when
$n>\max\{N_1,N_2\}$ and $0<\delta<\delta_1$, we have
\begin{equation*}
\big|\int\limits_{|x|\leq M}|x|dF_{n}(x)-\int\limits_{|x|\leq M}|x|dF_{sc}(x)\big|<\frac{\varepsilon}{3}+\frac{\varepsilon}{2}<\varepsilon
\end{equation*}
Thus, we have
\begin{equation}\label{eq1}
\lim\limits_{n\rightarrow \infty}\int\limits_{|x|\leq M}f(x)dF_{n}(x) =\int\limits_{|x|\leq M}f(x)dF_{sc}(x)\,\,\,\,a.s.
\end{equation}

By the similar method, we can prove that
\begin{equation*}
\lim\limits_{n\rightarrow \infty}\int\limits_{|x|\leq M}x^2 dF_{n}(x) =\int\limits_{|x|\leq M}x^2dF_{sc}(x)\,\,\,\,a.s.
\end{equation*}
Since
\begin{equation*}
\lim\limits_{n\rightarrow \infty}\int_{R}x^2 dF_{n}(x) =\int_{R}x^2dF_{sc}(x)\,\,\,\,a.s.
\end{equation*}
We obtain
\begin{equation*}
\lim\limits_{n\rightarrow \infty}\int\limits_{|x|>M}x^2 dF_{n}(x) =\int\limits_{|x|>M}x^2dF_{sc}(x)\,\,\,\,a.s.
\end{equation*}
Moreover, we know that $\int\limits_{|x|>M}x^2dF_{sc}(x)=0$ and $0\leq \int\limits_{|x|>M}|x| dF_{n}(x)\leq \int\limits_{|x|>M}x^2 dF_{n}(x)$.
It follows that
\begin{equation}\label{eq2}
\lim\limits_{n\rightarrow \infty}\int\limits_{|x|>M}|x| dF_{n}(x) =0=\int\limits_{|x|>M}|x|dF_{sc}(x)\,\,\,\,a.s.
\end{equation}
Combining Equations \ref{eq1} and \ref{eq2}, we have
\begin{equation}\label{eq3}
\lim\limits_{n\rightarrow \infty}\int_{R}|x| dF_{n}(x) =\int_{R}|x|dF_{sc}(x)\,\,\,\,a.s.
\end{equation}

Now, it is time to estimate the matching energy of random graphs.
\begin{eqnarray*}
\frac{ME(G_{n,p})}{n^{3/2}p^{1/2}}&=&\frac{1}{n^{3/2}p^{1/2}}\sum\limits_{i=1}^{n}|x_i|=\frac{1}{n}\sum\limits_{i=1}^{n}|\lambda_i|\\
&=&\int |x|\,d F_{n}(x)\stackrel{n\rightarrow \infty}{\longrightarrow}\int |x|\,\rho_{sc}(x)\,dx\quad a.s.\\
&=&\frac{1}{2\pi}\int_{-2}^{2} |x|\sqrt{4-x^2}dx\\
&=&\frac{8}{3\pi}.
\end{eqnarray*}
The proof is thus complete. \qed

\begin{appendix}
\section{Theorem \ref{THM0} implies Theorem \ref{THM1}}\label{Appendix}
\pf Let $f$ be a bounded continuous function on $R$. Let $M>4$ and $\delta>0$.
By the Weierstrass approximation theorem, there exists a polynomial $P_\delta$
such that $|f(x)-P_{\delta}(x)|<\delta$ for all $|x|\leq M$.
\begin{equation}\label{equa1}
\big|\int_R fdF_n(x)-\int_R fdF_{sc}(x)\big|\leq |\int_R P_{\delta}d(F_{n}(x)-F_{sc}(x))|
+2\delta+\big|\int\limits_{[-M,M]^c} (f-P_{\delta}) dF_{n}(x)\big|
\end{equation}
Suppose the degree of $P_{\delta}$ is $a$. Since $f$ is bounded, we have $f(x)<Cx^a$,
for $|x|\geq M$ and some constant $C$. Then we have
\begin{equation*}
\big|\int\limits_{[-M,M]^c} (f-P_{\delta}) dF_{n}(x)\big|\leq\big|\int\limits_{[-M,M]^c} C'|x^a|dF_{n}(x)\big|\leq C'\frac{1}{M^a}\int_R x^{2a}dF_{n}(x)
\end{equation*}
where $C'$ is a constant.

As $n\rightarrow \infty$, we have that the first part of Equation \ref{equa1} converges to zero, and the last part
\begin{equation*}
\lim\limits_{n\rightarrow \infty} C'\frac{1}{M^a}\int x^{2a}dF_{n}(x)\leq C'\frac{1}{M^a}\int_{-2}^{2}x^{2a}dF_{sc}(x)\leq  4C'\frac{4^a}{M^a}
\end{equation*}
By Theorem \ref{THM0}, for every $\epsilon>0$, there exists a number $N_1$ such that for $n\geq N_1$,
\begin{equation}\label{equa2}
|\int P_{\delta}d(F_{n}(x)-F_{sc}(x))|<\epsilon/3.
\end{equation}
There exists a number $M_1$, such that for $M>M_1$,
\begin{equation}\label{equa3}
4C'\frac{4^a}{M^a}<\epsilon/3.
\end{equation}
Let $\delta<\frac{\epsilon}{6}$. By Equations (\ref{equa1}), (\ref{equa2}) and (\ref{equa3}), we have
\begin{equation*}
\big|\int fdF_n(x)-\int fdF_{sc}(x)\big|<\epsilon\,\,.
\end{equation*}
The proof is thus complete. \qed
\end{appendix}
\end{document}